\documentclass[11pt]{article}

\usepackage{amsmath, amssymb, amscd}

\def\eqref#1{(\ref{#1})}
\newcommand{\goth}{\frak}

\newcommand{\arrow}{{\:\longrightarrow\:}}
\newcommand{\Z}{{\Bbb Z}}
\newcommand{\C}{{\Bbb C}}
\newcommand{\R}{{\Bbb R}}

\newcommand{\6}{\partial}
\def\1{\sqrt{-1}\:}

\newcommand{\cntrct}                
{\hspace{2pt}\raisebox{1pt}{\text{$\lrcorner$}}\hspace{2pt}}

\newcommand{\calo}{{\cal O}}

\renewcommand{\bar}{\overline}
\renewcommand{\phi}{\varphi}
\renewcommand{\epsilon}{\varepsilon}
\renewcommand{\geq}{\geqslant}
\renewcommand{\leq}{\leqslant}

\newcommand{\const}{{\it const}}

\newcommand{\End}{\operatorname{End}}
\newcommand{\Coh}{\operatorname{Coh}}

\newcommand{\Id}{\operatorname{Id}}

\newcommand{\Vol}{\operatorname{Vol}}

\newcommand{\Aut}{\operatorname{Aut}}
\newcommand{\Iso}{\operatorname{Iso}}

\newcommand{\codim}{\operatorname{codim}}

\newcommand{\slope}{\operatorname{slope}}
\newcommand{\rk}{\operatorname{rk}}

\newcommand{\Lie}{\operatorname{Lie}}

\newcommand{\Tr}{\operatorname{Tr}}

\newcommand{\comment}[1]{{}}

\def\blacksquare{\hbox{\vrule width 4pt height 4pt depth 0pt}}
\def\endproof{\blacksquare}

\makeatletter

\@ifundefined{Bbb}
     {\newcommand{\Bbb}[1]{{\mathbb #1}}}%
{}

\newcommand{\ps@verbit}{%
  \renewcommand{\@oddhead}{%
          \scriptsize
          {Holomorphic bundles on Hopf manifolds}
          \hfil\tiny {M. Verbitsky, August 27, 2004 }}
  \renewcommand{\@evenhead}{\@oddhead}
  \renewcommand{\@oddfoot}{\hfil\thepage\hfil}
  \renewcommand{\@evenfoot}{\@oddfoot}}
 
\newcounter{Mycounter}[section]
\newcounter{lemma}[section]
\renewcommand{\thelemma}{\noindent{Lemma \thesection.\arabic{lemma}}}
\newcommand{\lemma}{%
     \setcounter{lemma}{\value{Mycounter}}
     \refstepcounter{lemma}
     \stepcounter{Mycounter}
     {\bf \thelemma:\ }}

\newcounter{claim}[section]
\renewcommand{\theclaim}{\noindent{Claim \thesection.\arabic{claim}}}
\newcommand{\claim}{%
     \setcounter{claim}{\value{Mycounter}}
     \refstepcounter{claim}
     \stepcounter{Mycounter}
     {\bf \theclaim:\ }}

\newcounter{sublemma}[section]
\newcounter{corollary}[section]
\renewcommand{\thecorollary}{\noindent{Corollary \thesection.\arabic{corollary}}}
\newcommand{\corollary}{%
     \setcounter{corollary}{\value{Mycounter}}
     \refstepcounter{corollary}
     \stepcounter{Mycounter}
     {\bf \thecorollary:\ }}

\newcounter{theorem}[section]
\renewcommand{\thetheorem}{\noindent{Theorem \thesection.\arabic{theorem}}}
\newcommand{\theorem}{%
     \setcounter{theorem}{\value{Mycounter}}
     \refstepcounter{theorem}
     \stepcounter{Mycounter}
     {\bf \thetheorem:\ }}

\newcounter{conjecture}[section]
\newcounter{proposition}[section]
\renewcommand{\theproposition}
       {\noindent{Proposition \thesection.\arabic{proposition}}}
\newcommand{\proposition}{%
     \setcounter{proposition}{\value{Mycounter}}
     \refstepcounter{proposition}
     \stepcounter{Mycounter}
     {\bf \theproposition:\ }}

\newcounter{definition}[section]
\renewcommand{\thedefinition}
       {\noindent{Definition~\thesection.\arabic{definition}}}
\newcommand{\definition}{%
     \setcounter{definition}{\value{Mycounter}}
     \refstepcounter{definition}
     \stepcounter{Mycounter}
     {\bf \thedefinition:\ }}

\newcounter{example}[section]
\newcounter{remark}[section]
\renewcommand{\theremark}{\noindent{Remark \thesection.\arabic{remark}}}
\newcommand{\remark}{%
     \setcounter{remark}{\value{Mycounter}}
     \refstepcounter{remark}
     \stepcounter{Mycounter}
     {\bf \theremark:\ }}

\newcounter{problem}[section]
\newcounter{question}[section]
\@addtoreset{equation}{section}
\@addtoreset{footnote}{section}
\makeatother

\begin{document}

\begin{center}
{\LARGE\bf
Holomorphic bundles \\[3mm] on diagonal Hopf manifolds
}
\\[4mm]
Misha Verbitsky\footnote{Misha Verbitsky is 
an EPSRC advanced fellow 
supported by CRDF grant RM1-2354-MO02 
and EPSRC grant  GR/R77773/01}
\\[4mm]

{\tt verbit@maths.gla.ac.uk, \ \  verbit@mccme.ru}
\end{center}

{\small 
\hspace{0.15\linewidth}
\begin{minipage}[t]{0.7\linewidth}
{\bf Abstract} \\
Let $A\in GL(n, \C)$ be a diagonal linear operator,
with all eigenvalues satisfying $|\alpha_i|<1$,
and $M = (\C^n\backslash 0)/\langle A\rangle$
the corresponding Hopf manifold. We show that any
stable holomorphic bundle on $M$ can be lifted to 
a $\tilde G_F$-equivariant coherent sheaf on $\C^n$, 
where $\tilde G_F \cong (\C^*)^l$ is a commutative
Lie group acting on $\C^n$ and containing $A$.
This is used to show that all stable bundles
on $M$ are filtrable,
that is, admit a filtration by a sequence $F_i$ of 
coherent sheaves, with all subquotients 
$F_i/F_{i-1}$ of rank 1.
\end{minipage}
}

{
\small
\tableofcontents
}


\section{Introduction}
\label{_Intro_Section_}


In this paper we study the Hopf manifolds of form
$M = (\C^n\backslash 0)/\langle A\rangle$,
where $A\in GL(n, \C)$ is a linear operator
with all eigenvalues satisfying $|\alpha_i| <1$
(such an operator is called a linear contraction).
Deforming $A$ to an operator $\lambda \cdot Id$,
$0<|\lambda|<1$, we find that $M$ is diffeomorphic
to $S^{2n-1}\times (\R/\Z) \cong S^{2n-1}\times S^1$.
The odd Betti numbers of $M$ are odd, hence $M$ is not K\"ahler.
This is the first example of non-K\"ahler manifold 
known in algebraic geometry.

When $A$ is diagonal and has form $A = \tau \cdot Id$,
$M$ is elliptically fibered over $\C P^{n-1}$,
with all fibers isomorphic to an elliptic curve
$C_\tau = \C^*/\langle \tau\rangle$.
In this (so-called ``classical'') case, 
the algebraic dimension is maximal possible.
For arbitrary $A$, the algebraic dimension
of $M$ can reach any value from 0 to $n-1$.

Algebraic geometry of Hopf manifolds, especially
Hopf surfaces, is well studied (\cite{_Kato1_}, \cite{_Kato2_},
\cite{_Brinzanescu_Moraru:FM_}, \cite{_Brinzanescu_Moraru:stable_},
\cite{_Moraru:Hopf_}). For $\dim M =2$, one has 
a good understanding of the geometry of holomorphic
vector bundles on $M$ (\cite{_Moraru:Hopf2_}).
A typical stable vector bundle in this situation
is non-filtrable, and actually contains no proper
holomorphic subsheaves. 

For $\dim_\C M >2$, geometry of holomorphic
vector bundles is drastically different. 
In \cite{_Verbitsky:Sta_Elli_}, it was shown
that any bundle (and any coherent sheaf)
on a classical Hopf manifold 
\[ 
 (\C^n\backslash 0)/\langle \lambda \cdot Id\rangle, 
 \ \  n>2, \ \ 0< |\lambda|<1
\]
is filtrable. In the present paper, we generalize
this theorem to an arbitrary diagonal Hopf manifold.

\hfill

\theorem\label{_main_filtra_Theorem_}
Let $A\in GL(n, \C)$ be a diagonal
linear operator, with all eigenvalues
satisfying $|\alpha_i|<1$, and 
$M = (\C^n\backslash 0)/\langle A\rangle$
the corresponding Hopf manifold. Then any coherent
sheaf $F\in \Coh(M)$ is {\bf filtrable}, that is,
admits a filtration 
\[ 
0 = F_0 \subset F_1 \subset ... \subset F_m = F
\]
with $\rk F_i /F_{i-1} \leq 1$.

\hfill

{\bf Proof:} Using induction, we can always assume 
that any sheaf $F'$ with $\rk F' < \rk F$
is filtrable. Then $F$ is filtrable unless $F$
has no proper coherent shubsheaf. In the lattter case, $F$
is stable. Therefore, \ref{_main_filtra_Theorem_}
is implied by the following theorem, which is proven in
Section \ref{_equi_shea_Section_} by 
the means of gauge theory.

\hfill

\theorem\label{_stable_filtra_main_Theorem_}
Let $A\in GL(n, \C)$ be a diagonal
linear operator, with all eigenvalues
satisfying $|\alpha_i|<1$, and 
$M = (\C^n\backslash 0)/\langle A\rangle$
the corresponding Hopf manifold.
We choose a locally conformally K\"ahler
Hermitian structure on $M$
as in Subsection \ref{_LCK_Hopf_Subsection_}.
Let $F$ be a holomorphic bundle
(or a reflexive coherent sheaf)
which is stable with respect to this
Hermitian structure.\footnote{For a definition
of stability on Hermitian manifolds, see 
Section \ref{_stable_Kobaya_Section_}.}
Then $F$ is filtrable.

\hfill

{\bf Proof:} See \ref{_equi_shea_and_shea_on_Hopf_Remark_}.
\endproof

\hfill

The proof of \ref{_stable_filtra_main_Theorem_}
goes as follows. Using the 
Kobayashi-Hitchin correspondence 
on complex Hermitian manifolds
(Section \ref{_stable_Kobaya_Section_}),
we show that any stable bundle on a
diagonal Hopf manifold is equivariant
with respect to a certain holomorphic flow
(\ref{_all_stable_equi_Corollary_}).
Taking a completion of this flow in $GL(n, \C)$,
we obtain an abelian Lie group, which is
isomorphic to $(\C^*)^l$ (\ref{_tilde_G_C^*^l_Proposition_}). 
This allows us to treat stable holomorphic 
bundles (or reflexive sheaves) on $M$ as objects in a category
$(\C^*)^l$-equivariant coherent sheaves on $\C^n\backslash 0$
(\ref{_equi_shea_and_shea_on_Hopf_Remark_}).
Then we show that all objects in this category
are filtrable (\ref{_Coh_G_F_filtra_Theorem_}).


\section{Diagonal Hopf manifolds in Vaisman geometry}


\subsection{An introduction to Vaisman geometry}

\definition
Let $M$ be a complex manifold, $\dim_\C M >1$, and
$\tilde M$ its covering.
Assume that $\tilde M$ is equipped with a K\"ahler form
$\omega_K$, in such a way that the deck transform
of $\tilde M/M$ acts on $(\tilde M, \omega_K)$ 
by homotheties. 
The form $\omega_K$ defines on $M$ a conformal class by $[\omega_K]$.
The pair ($M,[\omega_K]$)  is called {\bf locally conformally
K\"ahler (LCK)}. A Hermitian form $\omega_H$ on $M$
is called an LCK-form if it belongs to the 
conformal class $[\omega_K]$.

\hfill

\definition
Consider an LCK-manifold $M$ with an LCK-form $\omega_H$.
A pullback of $\omega_H$ to $\tilde M$ 
is written as $f\omega_K$, where $f$
is a function and $\omega_K$ is the 
K\"ahler form on $\tilde M$. Therefore,
$d \omega_H = \omega_H \wedge \theta$,
where $\theta= \frac {d f} f$ is a 1-form on $M$. Clearly, 
$\theta$ is defined uniquely. Since $\theta=d\log f$, 
$\theta$ is also closed. This form is called 
{\bf the Lee form} of $(M, \omega_H)$.

\hfill

\remark
For a general Hermitian complex 
manifold $(M, \omega_H)$, the Lee form
is defined as ${d^c}^*\omega_H$, where $d^c =I \circ d
\circ I^{-1}$ is the twisted de Rham differential,
and ${d^c}^*$ its Hermitian adjoint. It is not difficult
to check that this definition is compatible with the 
one we used above. 

\hfill

\definition\label{_Gauduchin_defi_Definition_}
Let $(M, \omega_H)$ be a Hermitian complex manifold,
$\dim_\C M=n$.
Then $\omega_H$ is called {\bf a Gauduchon metric}
if $d^*{d^c}^*\omega_H=0$, or, equivalently,
$d d^c(\omega_H^{n-1}) =0$.

\hfill

\remark 
In \cite{_Gauduchon_1984_}, P. Gauduchon proved
that such a metric on $M$ exists and is unique,
up to a constant multiplier, in any
conformal class, provided that the
manifold $M$ is compact.

\hfill

\remark
On a compact LCK-manifold, this result translates
into an existence of a unique metric with a harmonic
Lee form $\theta$. Indeed, ${d^c}^*\omega_H=\theta$
is always closed, hence the Gauduchon condition
$d^*{d^c}^*\omega_H=0$ is equivalent to
$d^*\theta=0$. 

\hfill

Further on, we shall always fix a choice of a Hermitian
metric on an LCK-manifold by choosing a Gauduchon metric.

\hfill

\definition
Let $M$ be an LCK-manifold equipped with a Gauduchon
metric $\omega_H$, $\theta$ its Lee form and
$\nabla$ the Levi-Civita connection associated with
$\omega_H$. Assume that $\theta$ is parallel: 
$\nabla\theta=0$. Then $M$ is called {\bf a Vaisman
manifold.}

\hfill

\remark\label{_Kamishima_Ornea_Remark_}
According to  Kamishima-Ornea (\cite{_Kamishima_Ornea_}),
a compact LCK-manifold $M$ is Vaisman if and only if
it admits a holomorphic vector field acting on $M$
conformally, in such a way that its lifting
to $\tilde M$ is not an isometry of 
$(\tilde M, \omega_K)$.

\hfill

\remark 
It is easy to see (\cite{_Dragomir_Ornea_}) that
the condition $\nabla \theta=0$ implies that
the dual to $\theta$ vector field $\theta^\sharp$
(called {\bf the Lee field}) is a holomorphic isometry
of $M$ and acts on  $\tilde M$ by non-isometric
conformal automorphisms. This gives
the ``only if'' part of Kamishima-Ornea
theorem. 

\hfill

For further results, details and calculations in
Vaisman geometry, the reader is referred to 
\cite{_Dragomir_Ornea_}, \cite{_Gauduchon_Ornea_},
\cite{_OV:Structure_}, \cite{_OV:Immersion_}, \cite{_OV:Potential_}.

\hfill

Further on, we shall use the following lemma, which is
proven in \cite{_Verbitsky:LCHK_} 
(see also \cite{_OV:Immersion_}).

\hfill

\lemma \label{_omega_0_in_Vaisman_Lemma_}
Let $M$ be a Vaisman manifold, $\theta^\sharp$ it Lee
field, and $\Sigma$ the complex holomorphic foliation
generated by $\theta^\sharp$. Denote
by $\omega_0:= d^c \theta$ the real (1,1)-form
obtained as a $d^c = I \circ d \circ I^{-1}$-differential
of the Lee field $\theta$. Then $\omega_0\geq 0$,
and the null direction of $\omega_0$ is precisely $\Sigma$.

\endproof

\hfill

\remark \label{_weight_bu_Remark_}
Let $L_{\R}$ be a real flat line bundle on $M$ with the
same automorphy factors as the K\"ahler form $\omega_K$
(in conformal geometry, it is known as 
{\bf the weight bundle}).
Any non-degenerate positive section of $L_{\R}$
corresponds uniquely to a metric on $M$
conformally equivalent to $\omega_K$, and
the converse is also true. The Gauduchon
metric gives a rise to a section $\mu_G$ of $L_{\R}$.
Consider $L:= L_{\R}\otimes_{\R} \C$
as a holomorphic Hermitian line bundle, with a 
holomorphic structure induced from the
flat connection on $L=L_{\R}\otimes_{\R} \C$,
and Hermitian structure defined by 
$|\mu_G|=const$. Denote by $\nabla_C$ the
corresponding Chern connection. Then
$\omega_0$ is the curvature of $\nabla_C$
(\cite{_Verbitsky:LCHK_}, \cite{_OV:Immersion_}).

\subsection{LCK structure on diagonal Hopf manifolds}
\label{_LCK_Hopf_Subsection_}

The main examples of LCK and Vaisman geometries
are provided by the theory of Hopf manifolds.

\hfill

\definition
Let $A\in GL(n)$ be a linear transform, acting on $\C^n$
with all eigenvalues satisfying $|\alpha_i|<1$.
Denote by $\langle A \rangle\subset GL(n, \C)$ the cyclic group
generated by $A$. The quotient 
$(\C^n\backslash 0)/\langle A \rangle$
is called {\bf a linear Hopf manifold}.
If $A$ is diagonalizable, $(\C^n\backslash 0)/\langle A \rangle$
is called {\bf a diagonal Hopf manifold}.

\hfill

\remark
If one takes an arbitrary holomorphic contraction
$A$ instead of a linear contraction, one obtains
the general definition of a Hopf manifold (see
e.g. \cite{_Kato1_}, \cite{_Kato2_} for details). 

\hfill

\remark
Izu Vaisman, who introduced the subject and studied 
the Vaisman manifolds at great length (see 
\cite{_Vaisman:Dedicata_}, \cite{_Vaisman:Torino_}),
called them the generalized Hopf manifolds.
This name is not suitable because many
Hopf manifolds are not Vaisman.
For linear Hopf manifolds, 
$(\C^n\backslash 0)/\langle A\rangle$ is Vaisman
if and only if $A$ is diagonalizable
(see \cite{_OV:Potential_}).

\hfill

Let $A\in GL(n, \C)$ be a diagonal linear transform:
\begin{equation*}
\begin{bmatrix} \alpha_1 & 0 & \dots & 0\\
0 & \alpha_2 & \dots & 0\\
\vdots & \vdots & \ddots & \vdots\\
0 & 0 & \dots & \alpha_n,
\end{bmatrix}, \hfill \hfill |\alpha_i| <1
\end{equation*}
Consider the K\"ahler metric $\omega_K:= -\1\6\bar\6 \phi$ on
$\C^n \backslash 0$, defined using the K\"ahler potential 
$\phi:\; \C^n \backslash 0\arrow \R$. The $\phi$ is defined
via the formula
\begin{equation}\label{_potential_phi_formula_Equation_} 
\phi(t_1, ... , t_n) = \sum |t_i|^{\beta_i},
\end{equation}
where $\beta_i:= \log_{|\alpha_i|^{-1}} C$
are positive real numbers which satisfy 
$|\alpha_i|^{-\beta_i}= C$
for some fixed real constant $C>1$, 
chosen in such a way that all $\beta_i$ 
satisfy $|\beta_i| \geq 2$, and $t_i$
are complex coordinates.
By construction, $A^*\phi = C^{-1} \phi$. Indeed,
\[ A^*\phi(t_1, ..., t_n) = \phi(A(t_1, ..., t_n))=
   \sum |\alpha_i t_i|^{\beta_i}= C^{-1}\phi(t_1, ..., t_n).
\]
Therefore, $\omega_K:= -\1\6\bar\6 \phi$
is a K\"ahler form which satisfies $A^* \omega_K =C^{-1}\omega_K$.
This implies that the diagonal Hopf manifold
$(\C^n\backslash 0)/\langle A \rangle$ is LCK.
To see that it is Vaisman, we notice that the
holomorphic vector field $\log A$ acts on 
$(\C^n\backslash 0, \omega_K)$ conformally and
apply the Kamishima-Ornea theorem 
(\ref{_Kamishima_Ornea_Remark_}). 

We proceed with computing the Lee field for the
Gauduchon metric on $(\C^n\backslash 0)/\langle A \rangle$,
equipped with a conformal structure defined by the K\"ahler
form described above. 

Consider the action of the complex Lie group
$V(t)= e^{\C v}$ generated by the holomorphic vector field
$v:= \sum_i -t_i\log |\alpha_i| \frac{d}{dt_i}$. 
By construction, $V(\lambda)$ is a linear operator
which can be written as $\sum e^{|\alpha_i|\lambda} t_i$.
For $\lambda$ real, this operator multiplies $\phi$
by a constant $C^\lambda$ (this is proven in the same was as
one proves that $A(\phi)=C^{-1}\phi$), and for $\lambda$ purely imaginary,
$V(\lambda)$ preserves $\phi$ (this is clear).
Therefore, $v^c:= I(v)$ acts on $(\C^n\backslash 0, \omega_K)$
by holomorphic isometries. 

The corresponding moment map $\mu:\; \C^n\backslash 0\arrow \R$
is given by $d\mu = \omega_K(v^c, \cdot)$. The latter differential
form is written as 
\begin{equation}\label{_omega_K(v^c)_Equation_}
(d d^c \phi)\cntrct v^c =
\Lie_{v^c} d^c \phi - d(d^c\phi \cntrct v^c).
\end{equation}
The first term of the right hand side of 
\eqref{_omega_K(v^c)_Equation_} vanishes because
$v^c$ acts on $(\C^n\backslash 0)$ preserving $\phi$
and a complex structure. This gives
\[ \omega_K(v^c, \cdot)=(d d^c \phi)\cntrct v^c = -d(d^c\phi \cntrct v^c)
= d(d\phi\cntrct v)= \log C \cdot d\phi
\]
(the last equation holds because 
$d\phi\cntrct v= \Lie_v \phi = \log C\cdot \phi$).
We obtained that $\log(C)\phi$ is the moment map for $V(t)$ acting on
$(\C^n\backslash 0, \omega_K)$.

We obtained the following claim, which is well known
in many similar situations. 

\hfill

\claim\label{_moment_map_Claim_}
Let $A\subset GL(n)$ be a diagonal contraction
of $\C^n$, with all eigenvalues $\alpha_i$ satisfying
$|\alpha_i| <1$. Consider a K\"ahler metric
$\omega_K:= -\1 \6\bar\6 \phi$ on $\C^n\backslash 0$,
where the K\"ahler potential $\phi$ is defined
by the formula \eqref{_potential_phi_formula_Equation_},
and let $V(t)= e^{\C v}$ be the holomorphic flow generated
by $v:=\sum_i -t_i\log |\alpha_i| \frac{d}{dt_i}$. 
Let $v^c:= I(v)$ be the complex adjoint of $v$.
Then $e^{\R v^c}\subset V(t)$, 
preserves the K\"ahler structure on $\C^n\backslash 0$,
and the corresponding moment map is $\log (C)\phi$:
\[ d(\log (C) \phi) = \omega_k(v^c, \cdot).
\]

\endproof

\hfill

Consider the Hermitian form $\omega_H = \frac{\omega_K}\phi$
on $M = (\C^n\backslash 0)/\langle A \rangle$.
The corresponding Lee form $\theta$ is obtained via 
\[ d\omega_H = - \frac{\omega_K}{\phi^2}= - \omega_H \wedge \log d\phi,
\]
hence $\theta= \frac{d\phi}{\phi}$.
The dual under $\omega_H$ vector field (Lee field) is given by
$\theta^{\sharp}= v$, where $v=\log (C)\sum_i -t_i\log |\alpha_i| \frac{d}{dt_i}$.
This is clear because $v$ is dual to $\log (C)d\phi$
with respect to $\omega_K$ as \ref{_moment_map_Claim_}
implies, and $\omega_H = \frac{\omega_K}\phi$.

This gives the following Proposition.

\hfill

\proposition\label{_Lee_field_on_Hopf_Proposition_}
In assumptions of \ref{_moment_map_Claim_},
consider the Hermitian form $\omega_H = \frac{\omega_K}\phi$
on $M = (\C^n\backslash 0)/\langle A \rangle$. 
Then the corresponding Lee field is given
as 
\begin{equation}\label{_theta_sharp_via_alpha_i_Equation_}
\theta^{\sharp}=\log (C)\sum_i -t_i\log |\alpha_i| \frac{d}{dt_i}.
\end{equation}
Moreover, $\omega_H$ is Gauduchon.

\hfill

{\bf Proof:} The equation
\eqref{_theta_sharp_via_alpha_i_Equation_}
is proven above. To see that $\omega_H$ is Gauduchon,
it suffices to see that $|\theta^\sharp|_{\omega_H}$
is constant. Indeed, from the definition of $d^*$
it follows easily that 
\[ d^*\theta = \nabla_{\theta^\sharp}\theta^\sharp.\]
However, $\theta^\sharp$ is Killing, 
because $\Lie_{\theta^\sharp} \phi = C \phi$,
$\Lie_{\theta^\sharp} \omega_K = C \omega_K$,
and therefore 
\[ \Lie_{\theta^\sharp} \omega_H = C \omega_H - C \omega_H=0.
\]
By another definition of Killing fields, this means that
\[ (\nabla_X\theta^\sharp, Y)_{\omega_H} = 
-(\nabla_Y\theta^\sharp, X)_{\omega_H}
\]
for all vector fields $X, Y$. 
Taking $Y = \theta^\sharp$,
and applying $\Lie_X(\theta^\sharp, \theta^\sharp)_{\omega_H}=0$,
we obtain
\[ 0 = (\nabla_X\theta^\sharp, \theta^\sharp)_{\omega_H}=
   -(\nabla_{\theta^\sharp}\theta^\sharp, X)_{\omega_H}.
\]
As $X$ is arbitrary,
this implies $\nabla_{\theta^\sharp}\theta^\sharp=0$.
Therefore, \ref{_Lee_field_on_Hopf_Proposition_}
is implied by the equation $\omega_H(\theta^\sharp, \bar \theta^\sharp)=const$,
or, equivalently,
\begin{equation}\label{_omega_K_theta_sharp_const_phi_Equation_}
\omega_K(\theta^\sharp, \bar \theta^\sharp)=const\cdot \phi.
\end{equation}
Writing $\omega_K$ as 
\[ \omega_K = -\1 \6\bar\6 \phi = 
   \sum_i dt_i \wedge d \bar t_i |t_i|^{\beta_i-2} \frac{\beta_i^2}{4},
\]
and using $\theta^\sharp = \log (C)\sum_i -t_i\log |\alpha_i| \frac{d}{dt_i}$,
we obtain
\begin{equation}\label{omega_k_of_theta_explic_Equation_}
\omega_K(\theta^\sharp, \bar \theta^\sharp) = 
\log (C)^2\sum_i (\log |\alpha_i|)^2 |t_i|^{\beta_i}\frac{\beta_i^2}{4}.
\end{equation}
By definition, $e^{-\log |\alpha_i| \beta_i}=C$,
in other words, $\beta_i = -\frac{\log C}{\log \alpha_i}$.
Plugging this into \eqref{omega_k_of_theta_explic_Equation_},
we obtain
\[ \omega_K(\theta^\sharp, \bar \theta^\sharp)= 
   \sum_i |t_i|^{\beta_i}\frac{(\log C)^4}{4} = \frac {(\log C)^4} 4 \phi.
\]
This proves \ref{_Lee_field_on_Hopf_Proposition_}. \endproof

\section{Stable bundles on Hermitian manifolds}
\label{_stable_Kobaya_Section_}

\subsection{Gauduchon metrics and stability}

\definition
Let $M$ be a compact complex Hermitian manifold.
Choose a Gauduchon metric in the same conformal 
class.\footnote{A Hermitian metric 
on a complex manifold of dimension $n$ is called {\bf Gauduchon}
if $\6\bar\6(\omega^{n-1})=0$, where $\omega$ is its
Hermitian form (\ref{_Gauduchin_defi_Definition_}). On a compact manifold,
a Gauduchon metric exists in any conformal class, and is unique
up to a constant multiplier, see \cite{_Gauduchon_1984_}.} 
Consider a torsion-free
coherent sheaf $F$ on $M$. Denote by $\det F$ its determinant
bundle. Pick a Hermitian metric $\nu$ on $\det F$, and
let $\Theta$ be the curvature of the associated Chern
connection. We define the degree of $F$ as follows:
\[
\deg F := \int_M \Theta \wedge \omega^{\dim_\C M-1},
\]
where $\omega\in \Lambda^{1,1}(M)$ is
the Hermitian form of the Gauduchon metric.
This notion is independent from the choice of the Hermitian
structure $\nu$ in $F$. 
Indeed, if $\nu' = e^\psi \nu$, $\psi \in C^\infty(M)$,
then the associated curvature form is written as
$\Theta' = \Theta + \6\bar\6 \psi$, and
\[ \int_M \6\bar\6 \psi\wedge \omega^{\dim_\C M-1}=0
\]
because $\omega$ is Gauduchon.

\hfill

If $F$ is a Hermitian vector bundle, $\Theta_F$ its
curvature, and the metric $\nu$ is induced from $F$, then
$\Theta= \Tr_F\Theta_F$. In K\"ahler case this allows
one to relate the degree of a bundle with the first Chern class.
However, in non-K\"ahler case, the degree is not a 
topological invariant --- it depends fundamentally 
on the holomorphic geometry of $F$. 
Moreover, the degree is not discrete, as in the K\"ahler 
situation, but takes values in continuum.

Further on, we shall
see that one can in some cases construct 
a holomorphic structure of any given 
degree $\lambda\in \R$ on a fixed 
$C^\infty$-bundle. In our examples, such
holomorphic structures are constructed
on a topologically trivial line bundle 
over a Vaisman manifold (\ref{_arbi_degree_Remark_}).

\hfill

\definition
Let $F$ be a non-zero torsion-free coherent sheaf on $M$.
Then $\slope(F)$ is defined as 
\[
\slope(F) := \frac{\deg F}{\rk F}.
\]
The sheaf $F$ is called\\
{\small
\begin{tabular}[t]{ll}
{\bf stable} & if for all subsheaves $F'\subset F$, 
we have $\slope(F')< \slope(F)$\\
{\bf semistable} & if for all subsheaves $F'\subset F$, 
we have $\slope(F')\leq \slope(F)$\\
{\bf polystable} & if $F$ 
can be represented as a direct sum of stable \\& coherent 
sheaves with the same slope. 
\end{tabular}
}

\hfill

\remark
This definition is stability is ``good'' as most
standard properties of stable and semistable bundles 
hold in this situation as well. In particular, all line bundles 
are stable; all stable sheaves are simple; the Jordan-H\"older
and Harder-Narasimhan filtrations are well defined
and behave in the same way as they do in
the usual K\"ahler situation 
(\cite{_Lubke_Teleman:Book_}, \cite{_Bruasse:Harder_Nara_}).

However, not all bundles are {\bf filtrable},
that is, are obtained as successive extensions
by coherent sheaves of rank 1. There are non-filtrable
holomorphic vector bundles on most non-algebraic K3 surfaces.

\subsection{Kobayashi-Hitchin correspondence}
\label{_Koba_Hi_Subsection_}

The statement of Kobayashi-Hitchin correspondence 
(Donaldson-Uhlenbeck-Yau theorem) is translated
to the Hermitian situation verbatim, following 
Li and Yau (\cite{_Li_Yau_}).

\hfill

\definition
Let $B$ be a holomorphic Hermitian vector
bundle on a Hermitian manifold $M$, and 
$\Theta\in \Lambda^{1,1}(M)\otimes \End(B)$
the curvature of its Chern connection $\nabla$. Consider the
operator $\Lambda:\; \Lambda^{1,1}(M)\otimes \End(B)\arrow \End(B)$
which is a Hermitian adjoint to $b \arrow \omega\otimes b$,
$\omega$ being the Hermitian form on $M$.
The connection $\nabla$ is called {\bf Hermitian-Einstein}
(or {\bf Yang-Mills}) if $\Lambda \Theta = \const \cdot \Id_B$.

\hfill

\theorem 
(Kobayashi-Hitchin correspondence)
Let $B$ be a holomorphic vector bundle on a compact 
complex manifold equipped with a Gauduchon metric. Then
$B$ admits a Hermitian-Einstein connection $\nabla$
if and only if $B$ is polystable. Moreover, the
Hermitian-Einstein connection is unique.

\hfill

{\bf Proof:} See \cite{_Li_Yau_}, \cite{_Lubke_Teleman:Book_}, 
\cite{_Lubke_Teleman:Universal_}. \endproof


\section{Stable bundles on Vaisman manifolds}


Existence of the positive exact (1,1)-form
$\omega_0$, defined in \ref{_omega_0_in_Vaisman_Lemma_},
brings many consequences for algebraic geometry of
the Vaisman manifolds (see e.g. \cite{_Verbitsky:LCHK_} and 
\cite{_OV:Immersion_}). One of these is the structure
theorem for Hermitian-Einstein bundles of degree 0.

The following result was stated and proven 
as Theorem 4.3, \cite{_Verbitsky:Sta_Elli_}
for positive principal elliptic fibrations,
which admit a similar structure. These manifolds
are not always Vaisman (e.g. Calabi-Eckmann manifolds
are not Vaisman). However, the proof of this theorem
can be repeated almost verbatim in the 
Vaisman situation.

\hfill

\theorem\label{_Hermi_Einste_curva_equi_Theorem_}
Let $M$ be a compact Vaisman manifold, $\dim_\C M >2$,
and $B$ a stable
bundle of degree 0 on $M$. Denote by $\Sigma$ the 1-dimensional
complex holomorphic foliation generated by the Lee field
$\theta^\sharp$. Then $\Theta(v, \cdot)=0$ for any 
$v\in \Sigma$.
In particular, $B$ is equivariant with respect to the
complex Lie group $V(t)$ generated by $\theta^\sharp$, and
this equivariant structure is compatible with the connection.

\hfill

{\bf Proof:} Consider the map 
\[ \Lambda:\; \Lambda^{1,1}(M, \End(B)) \arrow \End(B)\]
defined in Subsection \ref{_Koba_Hi_Subsection_}.
By definition, $\Theta$ is {\bf primitive},
that is, satisfies $\Lambda\Theta=0$. Then 
\ref{_Hermi_Einste_curva_equi_Theorem_}
is implied by the following proposition.

\hfill

\proposition\label{_primi_form_equi_Proposition_}
Let $M$ be a compact Vaisman manifold, $\dim_\C M >2$,
$B$ a Hermitian bundle with connection, and
$\Theta \in \Lambda^{1,1}(M, \R) \otimes_\R {\goth u}(B)$
a closed skew-Hermitian real (1,1)-form. Assume that $\Theta$ 
is primitive, that is, $\Lambda\Theta=0$.
Then $\Theta(v, \cdot) =0$ for any  $v\in \Sigma$.

\hfill

{\bf Proof:} Rescaling the metric, 
we normalize the Lee form $\theta$ so that $|\theta|=1$.
Let $\theta$, $\theta_1, ... , \theta_{n-1}$ be an orthonormal
basis in $\Lambda^{1,0}(M)$, with $\theta\in \Sigma$, 
$\theta_i\in \Sigma^\bot$.
Consider the form $\omega_0$ (\ref{_omega_0_in_Vaisman_Lemma_}).
This form is exact, positive, and has $n-1$
strictly positive eigenvalues. Using the basis described above,
we can write
\begin{equation}\label{_omega,_omega_0_explicit_Equation_}
\omega_H= -\1 
\left(\theta\wedge\bar\theta +\sum_{i}\theta_i \wedge \bar\theta_i \right), 
\ \ \  \omega_0 = -\1 \left (\sum_{i}\theta_i \wedge \bar\theta_i\right )
\end{equation}
where $\omega_H$ is the Hermitian form of $M$ 
(see \cite{_Verbitsky:LCHK_}, Proposition 6.1).

In this basis, we can write $\Theta$ as
\begin{align}\label{_Theta_basis_Equation_}
\Theta &= \sum_{i\neq j}(\theta_i \wedge \bar\theta_j
         + \bar\theta_i \wedge \bar\theta_j) \otimes b_{ij} +
\sum_{i}(\theta_i \wedge \bar\theta_i) \otimes a_i \\
& + \sum_{i}(\theta \wedge \bar\theta_i
         + \bar\theta \wedge \bar\theta_i) \otimes b_{i}
  + \theta\wedge\bar\theta\otimes a,
\end{align}
with $b_{ij}$, $b_i$, $a_i$, $a\in {\goth u}(B)$
being skew-Hermitian endomorphisms of $B$. 

Let $\Xi:= \Tr(\Theta\wedge \Theta)$. 
This is a closed (2,2)-form on $M$.
Then \eqref{_Theta_basis_Equation_} implies
\[
 (\1)^n\Xi \wedge \omega_0^{n-2} = \Tr\left(-\sum b_i^2 +a \left(\sum a_i\right)
\right)
\]
On the other hand, $ \sum a_i + a = \Lambda\Theta=0$,
hence 
\[
 (\1)^n\Xi \wedge \omega_0^{n-2} = \Tr\left(-\sum b_i^2 -a^2 \right).
\]
Since $u \arrow \Tr(-u^2)$ is a positive definite form on
${\goth u}(B)$, the integral
\begin{equation}\label{_integral_Xi_Equation_}
\int_M  (\1)^n\Xi \wedge \omega_0^{n-2}
\end{equation}
is non-negative, and positive unless $b_i$ and $a$ 
both vanish everywhere. 
Using $n>2$, we find that
\eqref{_integral_Xi_Equation_} vanishes, because
$\omega_0$ is exact and $\Xi$ is closed. Therefore,
$b_i$ and $a$ are identically zero, which is exactly the
claim of \ref{_primi_form_equi_Proposition_}. 
We proved \ref{_Hermi_Einste_curva_equi_Theorem_}.
\endproof

\hfill

\remark\label{_arbi_degree_Remark_}
The results of \ref{_Hermi_Einste_curva_equi_Theorem_}
can be applied to arbitrary stable bundle on $M$
using the following trick. Consider the line bundle
$L$ (\ref{_weight_bu_Remark_}). Write the Chern 
connection on $L$ as
\[ \nabla_C = \nabla_{triv} -\1\theta^c,\]
where $\theta^c= I(\theta)$ is the complex conjugate of $\theta$
(see \cite{_Verbitsky:LCHK_}, (6.11)), and $\nabla_{triv}$
is a trivial connection associated to the trivialization
of $L$ constructed in \ref{_weight_bu_Remark_}.
Since $d\theta^c = \omega_0$, 
$L$ has a degree $\delta:= \int \omega_0 \wedge \omega_H^{n-1}$
which is clearly positive (see 
\eqref{_omega,_omega_0_explicit_Equation_}). 
Given $\lambda\in \R$, 
denote by $L_\lambda$ a holomorphic Hermitian bundle
with the connection $\nabla_{triv} -\1\frac{\lambda}{\delta}\theta^c$.
Then $L_\lambda$ has degree $\lambda$.
We obtain that a Vaisman manifold admits
a line bundle $L_\lambda$ of arbitrary degree
$\lambda$. Moreover, $L_\lambda$ is by construction
$V(t)$-equivariant (the form $\theta^c$ is $V(t)$-invariant,
as $V(t)$ acts on $M$ preserving the metric and the
holomorphic structure). This brings the following corollary.

\hfill

\corollary\label{_all_stable_equi_Corollary_}
Let $M$ be a compact Vaisman manifold, and $B$ a stable
bundle. Consider a complex holomorphic flow $V(t)= e^{t\theta^\sharp}$
generated by the Lee field $\theta^\sharp$. Then $B$ admits
a natural $V(t)$-equivariant structure. 

\hfill

{\bf Proof:} Tensoring $B$ by $L_{\lambda}$ for appropriate
choice of $\lambda\in \R$, we obtain a stable bundle
of degree 0. Then \ref{_Hermi_Einste_curva_equi_Theorem_}
implies \ref{_all_stable_equi_Corollary_}. \endproof

\hfill


\section{Stable bundles on  Hopf 
manifolds and coherent sheaves on $\C^n$}


\subsection{Admissible Hermitian structures on reflexive sheaves}

\definition\label{_refle_Definition_}
Let $X$ be a complex manifold, and $F$ a coherent sheaf on $X$.
Consider the sheaf $F^*:= {\cal H}om_{\calo_X}(F, \calo_X)$.
There is a natural functorial map 
$\rho_F:\; F \arrow F^{**}$. The sheaf $F^{**}$
is called {\bf a reflexive hull}, or {\bf 
reflexization},
of $F$. The sheaf $F$ is called {\bf reflexive} if the map
$\rho_F:\; F \arrow F^{**}$ is an isomorphism. 

\hfill

\remark
For all coherent sheaves $F$, the map
$\rho_{F^*}:\; F^* \arrow F^{***}$ is an isomorphism
(\cite{_OSS_}, Ch. II, the proof of Lemma 1.1.12).
Therefore, a 
reflexive hull of a sheaf is always 
reflexive.

\hfill

Reflexive hull can be obtained by 
restricting to an open subset and taking the
pushforward.

\hfill

\lemma\label{_refle_pushfor_Lemma_}
Let $X$ be a complex manifold, $F$ a coherent sheaf on $X$,
$Z$ a closed analytic subvariety, $\codim Z\geq 2$, and
$j:\; (X\backslash Z) \hookrightarrow X$ the natural
embedding. Assume that the pullback $j^* F$ is
reflexive on $(X\backslash Z)$. Then the pushforward
$j_* j^* F$ is also reflexive. 

\hfill

{\bf Proof:} This is \cite{_OSS_}, Ch. II, Lemma 1.1.12.
\endproof

\hfill

\remark\label{_refle_pushfor_Remark_}
From \ref{_refle_pushfor_Lemma_},
it is apparent that one could obtain a reflexization
of a non-singular in codimension 1 coherent sheaf $F$
by taking $j_* j^* F$, where 
$j:\; (X\backslash Z) \hookrightarrow X$ the natural
open embedding, and $Z$ the singular locus of $F$.

\hfill

Using the results of \cite{_Bando_Siu_}, we are able
to apply the Kobayashi-Hitchin correspondence to
reflexive sheaves.

\hfill

\definition\label{_admissi_Definition_}
\cite{_Bando_Siu_} 
Let $F$ be a coherent sheaf on $M$ and
$\nabla$ a Hermitian connection
on $F$ defined outside of its singularities.
Denote by $\Theta$ the curvature of $\nabla$.
Then $\nabla$ is called {\bf admissible}
if the following holds
\begin{description}
\item[(i)] $\Lambda \Theta\in \End(F)$ is uniformly bounded
\item[(ii)] $|\Theta|^2$ is integrable on $M$. 
\end{description}

\hfill

\theorem \label{_B_S_exie_admissi_Theorem_}
\cite{_Bando_Siu_} Any torsion-free 
coherent sheaf admits an  admissible connection.
An admissible connection can be extended over the
place where $F$ is smooth. Moreover, if a bundle $B$ 
on $M\backslash Z$, $\codim_\C Z\geq 2$ is equipped
with an admissible connection, then $B$ can be extended to
a coherent sheaf on $M$. \endproof

\hfill

A version of Donaldson-Uhlenbeck-Yau theorem exists for 
coherent she\-a\-ves (\ref{_UY_for_shea_Theorem_});
given a torsion-free coherent sheaf $F$, $F$ admits an admissible
Hermitian-Einstein connection $\nabla$ if and only if $F$ is polystable.

\hfill

\theorem\label{_UY_for_shea_Theorem_}
Let $M$ be a compact K\"ahler manifold,  and $F$ a coherent
sheaf without torsion. Then $F$ admits an admissible 
 Hermitian-Einstein
metric is and only if $F$ is polystable. Moreover, if $F$
is stable, then this metric is unique, up to a constant
multiplier.

{\bf Proof:} \cite{_Bando_Siu_}, Theorem 3.
\endproof

\hfill

This proof can be adapted for Hermitian complex
manifolds with Gauduchon metric.

\subsection{Hermitian-Einstein bundles on Hopf manifolds and admissibility}

\theorem \label{_from_stable_to_refle_Theorem_}
Let $M= (\C^n\backslash 0)/\langle A\rangle$ be a diagonal
Hopf manifold, $n\geq 3$, and $B$ a stable holomorphic
bundle on $M$ of degree 0. Denote by $\tilde B$ the pullback of $B$ to
$\C^n \backslash 0$. Then $\tilde B$ can be extended
to a reflexive coherent sheaf $F$ on $\C^n$. Moreover, $F$ 
is $V(t)$-equivariant, where $V(t)$,
is the complex holomorphic flow on $\C^n$ generated 
by the Lee field $\theta^\sharp=\log (C)\sum_i -t_i\log |\alpha_i| \frac{d}{dt_i}$.

\hfill

{\bf Proof:} Consider a Hermitian-Einstein metric on $B$,
and lift it to $\tilde B$. Denote by $\tilde \Theta$ the
curvature of $\tilde B$. To extend $\tilde B$ to $\C^n$, we apply
the Bando-Siu theorem (\ref{_B_S_exie_admissi_Theorem_}). 
We need to show that $\tilde B$
is admissible, in the sense of \ref{_admissi_Definition_}. 
The K\"ahler metric
$\omega_K$ on $\C^n$ is conformally equivalent to that lifted from $M$,
hence $\Lambda \tilde \Theta=0$ (this condition means
that $\tilde \Theta$ is orthogonal to the Hermitian form 
pointwise, and therefore it is conformally invariant). 
To prove that $\tilde B$ is admissible, it
remains to show that $\tilde \Theta$ is square-integrable.
The function $|\tilde \Theta|^2$ can be expressed, using the
Hodge-Riemann relations, as follows. 

\hfill

\lemma \label{_tilde_Theta_Lemma_}
Let $B_1$ be a Hermitian bundle on a Hermitian
almost complex manifold $M_1$, of dimension $n$, and 
\[ \nu \in \Lambda^{1,1}(M_1, \goth{su}(B_1))
\]
a $\goth{su}(B_1)$-valued (1,1)-form satisfying $\Lambda(\nu)=0$.
Then 
\begin{equation}\label{_square_form_via_H_E_Equation_} 
|\nu|^2 = -\1\frac{n-1}{2n} \Tr (\Lambda^2 (\nu \wedge \nu)),
\end{equation}
where 
\[
  \Lambda:\; \Lambda^{p,q}(M_1, \goth{su}(B_1))
  \arrow \Lambda^{p-1,q-1}(M_1, \goth{su}(B_1))
\]
is the standard Hodge operator on differential 
forms. 

\hfill

{\bf Proof:} An elementary calculation,
and essentally the same as one which proves
the Hodge-Riemann bilinear relations (see e.g. \cite{_Bando_Siu_}).
\endproof

\hfill

\remark 
The equation \eqref{_square_form_via_H_E_Equation_} 
can be stated as 
\begin{equation}\label{_square_form_volume_via_H_E_Equation_} 
|\nu|^2\Vol(M_1) = 
-\1\frac{n-1}{2n\cdot 2^n\cdot n!} \Tr (\nu \wedge \nu)\wedge \omega_1^{n-2},
\end{equation}
where $\omega$ is the Hermitian form on $M_1$,
and $\Vol(M_1)$ the Riemannian volume.
This is clear from the definition of $\Lambda$ and
the relation $\Vol(M_1) = \frac1 {2^n n!} \omega_1^n$.

\hfill

Using \eqref{_square_form_volume_via_H_E_Equation_},
we obtain that $L^2$-integrability of $\tilde \Theta$
is equivalent to integrability of the form 
\begin{equation}\label{_square_curv_Equation_}
\Tr (\tilde \Theta\wedge \tilde \Theta)\wedge \omega_K^{n-2}.
\end{equation}

\hfill

The form $\tilde \Theta$ is by construction
$A$-invariant, and $\omega_K$ satisfies
$A^*(\omega_K) = c \omega_K$ because $M$ is LCK.
Therefore, the form \eqref{_square_curv_Equation_}
is homogeneous with respect to the action of $A$:
\begin{equation}\label{_curv_volu_homoge_Equation_}
A^* \left(\Tr (\tilde \Theta\wedge \tilde \Theta)\wedge \omega_K^{n-2}\right)=
    c^{n-2}\Tr (\tilde \Theta\wedge \tilde \Theta)\wedge \omega_K^{n-2},
  c<1. 
\end{equation}
Denote by $D$ the fundamental domain
for $\langle A\rangle$,
\[ D:= \{ x\in \C^n\backslash 0\ | \  1\leq \rho(x) <C\}
\]
Then $\C^n \backslash 0 = \cap_{i\in \Z} A^i(D)$.
To check that $\tilde \Theta$ is $L^2$-integrable
in a neighbourhood of 0, we need to show that the series
\[
\sum_{i=0}^{\infty} \int_{A^i(D)} |\tilde \Theta|^2 \Vol =
-\1\frac{n-1}{2n\cdot 2^n \cdot n!}\sum_{i=0}^{\infty} \int_{A^i(D)} 
 \Tr (\tilde \Theta\wedge \tilde \Theta)\wedge \omega_K^{n-2}
\]
converges. However, by homogeneity, the latter integral is
power series, and \eqref{_curv_volu_homoge_Equation_}
implies that it converges whenever $n>2$.
We have shown that $\tilde B$ is admissible.
Now, Bando-Siu theorem (\ref{_B_S_exie_admissi_Theorem_}) 
implies the first assertion of \ref{_from_stable_to_refle_Theorem_}.
The second assertion is implied immediately by \ref{_refle_pushfor_Lemma_}.
Indeed, let $(\C^n \backslash 0)\stackrel j \hookrightarrow \C^n$
be the standard embedding. Then $F = j_* \tilde B$
(\ref{_refle_pushfor_Lemma_}). By 
\ref{_all_stable_equi_Corollary_}, 
$\tilde B$ is $V(t)$-equivariant.
Then $j_* \tilde B$ is also $V(t)$-equivariant. 
\endproof 

\hfill

\remark\label{_exte_equiva_theorem_to_refle_Remark_}
Using the Bando-Siu version of Donaldson-Uhlenbeck-Yau theorem,
we can extend \ref{_from_stable_to_refle_Theorem_} verbatim
to reflexive coherent sheaves. 


\section{Equivariant sheaves on $\C^n$}
\label{_equi_shea_Section_}


\subsection{Extending $V(t)$-equivariance to 
$(\C^*)^l$-equivariance}

Let $M= (\C^n\backslash 0)/\langle A\rangle$ be a diagonal
Hopf manifold, and $V(t)= e^{\C \theta^\sharp}$, $t\in \C$ 
the holomorphic
flow generated by the Lee field $\theta^\sharp$ as above. 
Then $V(t)$ acts on $M$ by 
holomorphic isometries
(\cite{_Kamishima_Ornea_}). Consider the
closure $G$ of $V(t)$, $t\in \C$, within the group
$\Iso(M)$ of isometries of $M$. Denote by $\tilde G$
the the lifting of $G$ to $\Aut(\tilde M)$
(\cite{_OV:Structure_}, \cite{_OV:Immersion_}).
By construction, $\tilde G$ is the smallest 
closed Lie subgroup of $GL(n, \C)$ containing 
$V(t)$ and $A$. It is easy to check that $\tilde G$ is a reductive 
complex commutative Lie group. A similar result is
true for all Vaisman manifolds. 

\hfill

\proposition\label{_tilde_G_C^*^l_Proposition_}
For any Vaisman manifold $M$, let $\theta^\sharp$
be its Lee field, $G$ the closure of the corresponding
complex holomorphic 
flow within $\Iso(M)$, and $\tilde G$ its lift to
$\Aut(\tilde M)$. Then $\tilde G\cong (\C^*)^k$,
and the deck transform map $\gamma\in \Aut(\tilde M, M)$
lies in $\tilde G$. 

\hfill

{\bf Proof:} This is \cite{_OV:Immersion_}, 
Proposition 4.3.\endproof

\hfill

\proposition\label{_tilde_G_F_equiv_shea_Theorem_}
In assumptions of \ref{_from_stable_to_refle_Theorem_},
consider the action of the group $V(t)$ on 
$\Gamma(\C^n, F)$. Consider the adic topology
on $\calo_{\C^n}$ and $\Gamma(\C^n, F)$,
with $\lim f_i \arrow 0$ as $[f_i]_0\arrow \infty$,
where $[f_i]_0$ denotes the order of zeroes of
$f_i$ in $0\in \C^n$. Clearly, $V(t)$ is
continuous in adic topology. Let $\tilde G_F$
be the closure of $V(t)$-action on
$\Gamma(\C^n, F) \times \calo_{\C^n}$
in adic topology. Then
\begin{description}
\item[(i)] The natural map
$\tilde G_F \stackrel \rho \arrow GL(F/\goth m F)\times GL(\goth m/\goth m^2)$
is injective, where $\goth m$ is the maximal ideal of $0$ in
$\calo_{\C^n}$.

\item[(ii)] $\tilde G_F$ is a closure of
$V(t)$ under the natural map
$V(t)\arrow GL(F/\goth m F)\times GL(\goth m/\goth m^2)$.

\item[(iii)] Consider the natural projection
$\tilde G_F\stackrel \pi \arrow \tilde G$ induced
by 
\[ GL(F/\goth m F)\times GL(\goth m/\goth m^2)\arrow  GL(\goth m/\goth m^2).
\]
Then $\pi$ satisfies $g(af) = \pi(g)(a)g(f)$,
for any $f\in \Gamma(\C^n, F)$, $a \in \calo_{\C^n}$,
$g\in \tilde G_F$.
This gives a $\tilde G_F$-equivariant structure on $F$.

\item[(iv)] The group $\tilde G_F$ is isomorphic to 
$(\C^*)^l$.

\end{description}

{\bf Proof:} \ref{_tilde_G_F_equiv_shea_Theorem_} (i)
is clear from Nakayama's lemma. 
\ref{_tilde_G_F_equiv_shea_Theorem_} (ii)
is immediately implied by \ref{_tilde_G_F_equiv_shea_Theorem_} (i).
\ref{_tilde_G_F_equiv_shea_Theorem_} (iii) follows from
\ref{_tilde_G_F_equiv_shea_Theorem_} (ii) and $V(t)$-equivariance
of $F$. 

To prove \ref{_tilde_G_F_equiv_shea_Theorem_} (iv), we use
\ref{_tilde_G_F_equiv_shea_Theorem_} (ii), and notice that
$\tilde G_F$ is commutative as a closure of a 1-parametric
group within a Lie group $GL(F/{\goth m}F)\times GL(\goth m/\goth m^2)$. 
To show that 
$\tilde G_F\cong (\C^*)^l$, we need to prove that
it is reductive, that is, to show that $V(t)$
acts diagonally on $(F/\goth m F)\times(\goth m/\goth m^2)$.

The group $V(t)$ acts on $M$ holomorphically and conformally.
 Since the Hermitian-Einstein metric on $B$ is unique,  
up to a constant multiplier, the group $V(t)$ acts on 
$B$ also conformally. Then, $V(t)$ acts conformally
on the Hermitian space $\Gamma(B_{\C^n}, F)$
of holomorphic sections of $F$ on an open 
ball $B_{\C^n}\subset\C^n$ and on $\Gamma(B_{\C^n}, \calo_{\C^n})$. 
Since orthogonal matrices in finite dimension are diagonalizable,
  $V(t)$ acts diagonally on any finite-dimensional
subspace in $\Gamma(B_{\C^n}, F)\times \Gamma(B_{\C^n}, \calo_{\C^n})$
preserved by $V(t)$. Using the same classical
Poincare-Dulac argument as used in the proof
of Theorem 3.3 in \cite{_OV:Potential_}, we find
that $\Gamma(B_{\C^n}, F)\times \Gamma(B_{\C^n}, \calo_{\C^n})$
contains a dense (in appropriate, e.g. $\goth m$-adic topology) subspace
which is generated by finite-dimensional
$V(t)$-invariant subspaces. Then $V(t)$-action on the space
$\Gamma(B_{\C^n}, F)\times \Gamma(B_{\C^n}, \calo_{\C^n})$
is diagonal in a dense subspace. Therefore, 
this action is diagonal on its quotient
$(F/\goth m F)\times(\goth m/\goth m^2)$.
We proved \ref{_tilde_G_F_equiv_shea_Theorem_} (iv).
\endproof

\hfill

\remark\label{_gene_equi_to_shea_Remark_}
Using the Bando-Siu version of Donaldson-Uhlenbeck-Yau theorem
(see \ref{_exte_equiva_theorem_to_refle_Remark_}),
we can extend \ref{_tilde_G_F_equiv_shea_Theorem_} verbatim
to reflexive coherent sheaves.

\hfill

\remark\label{_equi_shea_and_shea_on_Hopf_Remark_}
Denote by $\C^n_*$ the complex manifold $\C^n\backslash 0$. 
Given a $\tilde G_F$-equivariant coherent sheaf on $\C^n_*$,
we can obtain a coherent sheaf on $\C^n_*/\langle A\rangle$.
Indeed, coherent sheaves on  $\C^n_* /\langle A\rangle$
are the same as $\langle A\rangle$-equivariant sheaves
on $\C^n_*$, and $\langle A\rangle$ lies in $\tilde G$
as \ref{_tilde_G_C^*^l_Proposition_} implies.
Therefore, to prove the filtrability of 
a stable bundle $B$ on $M = (\C^n_*) /\langle A\rangle,$
it suffices to show that the corresponding
$\tilde G_F$-equivariant coherent sheaf
$F$ is filtrable on $\C^n_*$  in the category
$\Coh_{\tilde G_F}(\C^n_*)$ of $\tilde G_F$-equivariant 
coherent sheaves. Then, the following theorem
proves \ref{_stable_filtra_main_Theorem_}.

\hfill

\theorem \label{_Coh_G_F_filtra_Theorem_}
Let $\tilde G_F\cong (\C^*)^l$ be a commutative Lie group,
acting on $\C^n_*$ via a homomorphism $\tilde G_F\stackrel \pi \arrow GL(\C, n)$,
and $\Coh_{\tilde G_F}(\C^n_*)$ be the category of
$\tilde G_F$-equivariant coherent sheaves on $\C^n_*$.
Assume that $\pi(\tilde G_F)$ contains an endomorphism
with all eigenvalues $<1$.
Then all objects of $\Coh_{\tilde G_F}(\C^n_*)$
are filtrable by $\tilde G_F$-equivariant coherent sheaves
of rank at most 1.

\hfill

We prove \ref {_Coh_G_F_filtra_Theorem_}
in Subsection \ref{_C^*^l_equiv_Subsection_}.

\subsection{$(\C^*)^l$-equivariant coherent sheaves on $\C^n\backslash 0$}
\label{_C^*^l_equiv_Subsection_}

We work in assumptions of \ref{_Coh_G_F_filtra_Theorem_}.

\hfill

\lemma\label{_R_gene_by_fini_G_F_inv_Lemma_}
Let $R\in \Coh_{\tilde G_F}(\C^n_*)$
be a $\tilde G_F$-equivariant coherent sheaf
over $\C^n_*:= \C^n\backslash 0$.
Then $R$ is generated over 
$\calo_{\C^n_*}$ by a finite-dimensional $\tilde G_F$-invariant 
space $V\subset \Gamma(R, \C^n_*)$.

\hfill

{\bf Proof:} The images of $\C^*$ are dense in
$\tilde G_F\cong (\C^*)^l$. Therefore, there
exists an embedding $\C^*\stackrel \mu \hookrightarrow \tilde G_F$
acting on $\C^n$ with all eigenvalues different from 1. 
This action can be written as 
\begin{equation*}
t \arrow \begin{bmatrix} t^{k_1} & 0 & \dots & 0\\
0 & t^{k_2} & \dots & 0\\
\vdots & \vdots & \ddots & \vdots\\
0 & 0 & \dots & t^{k_n}
\end{bmatrix}
\end{equation*}
where all $k_i$ are integers 
different from 0.
Clearly, $\mu$ acts on $\C^n_*$ freely in generic point, and the
quotient $\C^n_* /\mu (\C^*)$
is well defined. This quotient is known as 
a {\bf weighted projective space}, denoted by 
$\C P^{n-1}(k_1, k_2, ... k_n)$,
and it is a projective orbifold.
To give a $\mu$-equivariant coherent sheaf on $\C^n_*$
is by definition the same as to give a coherent 
sheaf on  the orbifold $\C P^{n-1}(k_1, k_2, ... k_n)$.
Let $R_0$ be the sheaf on $\C P^{n-1}(k_1, k_2, ... k_n)$
corresponding to $R$, considered as a 
$\mu$-equivariant sheaf on $\C^n_*$.
The sections of $R_0\otimes \calo(i)$
correspond to the sections of $R$ on which
$\mu(\C^*)$ acts with the weight $i$. 
We obtain a sequence of finite-dimensional
subspaces 
\[ \Gamma(R_0\otimes \calo(i))\subset \Gamma(R).
\]
Since $\calo(1)$ is ample,
the sheaf $R_0\otimes \calo(N)$
is globally generated for $N$ sufficiently big
(here we use the Kodaira-Nakano theorem for orbifolds,
\cite{_Baily_}). Then $\oplus_{i\leq N} \Gamma(R_0\otimes \calo(i))$
will generate $\Gamma(R)$ over $\calo_{\C^n_*}$.
Since $\tilde G_F$ commutes with $\mu(\C^*)$,
the space $\Gamma(R_0\otimes \calo(i))\subset \Gamma(R)$
is $\tilde G_F$-invariant. 
This proves \ref{_R_gene_by_fini_G_F_inv_Lemma_}.
\endproof

\hfill

Now we can prove the filtrability of 
arbitrary $R\in \Coh_{\tilde G_F}(\C^n_*)$.
By \ref{_R_gene_by_fini_G_F_inv_Lemma_},
for any $R\in \Coh_{\tilde G_F}(\C^n_*)$,
there exists a surjective $\tilde G_F$-equivariant map
$R_1 \arrow R\arrow 0$, where $R_1 = \calo_{\C^n_*}\otimes_\C W$,
and $W$ is a finite-dimensional representation of $\tilde G_F$.
Since $\tilde G_F$ is commutative, $W= \oplus W_i$, where
$W_i$ are $\tilde G_F$-invariant 1-dimensional subspaces
of $W$. This gives an epimorphism
\[
\oplus (\calo_{\C^n_*}\otimes W_i) \arrow R
\]
where all the summands $\calo_{\C^n_*}\otimes W_i$
are $\tilde G_F$-equivariant line bundles. 
Then, $R$ is clearly filtrable within 
$\Coh_{\tilde G_F}(\C^n_*)$. This proves
\ref{_Coh_G_F_filtra_Theorem_}. 
\ref{_stable_filtra_main_Theorem_} is also proven. 
\endproof

\hfill

{\bf Acknowledgements:} I am grateful to Ruxandra Moraru
who posed the problem, and Dmitry Novikov for a private lecture
on normal forms and the Poincare-Dulac theorem.

{\scriptsize

\hfill
}
{
\small

\noindent {\sc Misha Verbitsky\\
University of Glasgow, Department of Mathematics, \\
15 University Gardens, Glasgow G12 8QW, Scotland.}\\
\ \\
{\sc  Institute of Theoretical and
Experimental Physics \\
B. Cheremushkinskaya, 25, Moscow, 117259, Russia }\\
\ \\
\tt verbit@maths.gla.ac.uk, \ \  verbit@mccme.ru 
}

\end{document}